# STABILIZATION BY MEANS OF APPROXIMATE PREDICTORS FOR SYSTEMS WITH DELAYED INPUT

**Iasson Karafyllis**

**Department of Environmental Engineering, Technical University of Crete,
73100, Chania, Greece
email: ikarafyl@enveng.tuc.gr**

**Abstract**

Sufficient conditions for global stabilization of nonlinear systems with delayed input by means of approximate predictors are presented. An approximate predictor is a mapping which approximates the exact values of the stabilizing input for the corresponding system with no delay. A systematic procedure for the construction of approximate predictors is provided for globally Lipschitz systems. The resulting stabilizing feedback can be implemented by means of a dynamic distributed delay feedback law. Illustrating examples show the efficiency of the proposed control strategy.

## 1. Introduction

In the present work, we consider the stabilization problem for nonlinear systems with input delays and measurement delays of the form:

$$\dot{x}(t) = f(x(t), u(t-r))$$
$$x(t) = (x_1(t), \ldots, x_n(t))' \in \Re^n, \ u(t) \in U \tag{1.1}$$

where $r \geq 0$ is a constant and $U \subseteq \Re^m$ is a closed convex set with $0 \in U$. More specifically, we want to address the feedback design problem for system (1.1) based on the knowledge of a feedback stabilizer $u = k(x)$ for system (1.1) with no delay, i.e. (1.1) with $r = 0$, or

$$\dot{x}(t) = f(x(t), u(t))$$
$$x(t) = (x_1(t), \ldots, x_n(t))' \in \Re^n, \ u(t) \in U \tag{1.2}$$

In the literature there are different ways of attacking this problem:

- One way is to apply the feedback law $u(t) = k(x(t))$ and obtain conditions which guarantee stability for the closed-loop system. Such conditions can be obtained by using Lyapunov, Razumikhin or small-gain arguments as in [8,11].
- A second way is to modify the feedback law $u(t) = k(x(t))$ by applying a predictor, i.e., a mapping $p(t)$ which guarantees that $p(t-r) = x(t)$. The predictor feedback $u(t) = k(p(t))$ will guarantee that $u(t-r) = k(x(t))$ and thus we will obtain the stability properties of the closed-loop system (1.2) with $u = k(x)$. This idea is classical in linear systems (the Smith predictor, see [9] and the references in [5]) and was extended recently in nonlinear systems in [5] by M. Krstic,
- Finally, another way is to exploit certain characteristics of the system in order to obtain a modified feedback law $\tilde{k}(x)$ such that the application of the modified feedback law $u(t) = \tilde{k}(x(t))$ will ensure stability for the corresponding closed-loop system (see [2,6,7]).

In the present work, we apply the "predictor approach" and we obtain results which guarantee global asymptotic stability for the closed-loop system for arbitrary large values of the delay $r$. The proofs of the main results of the present work rely heavily on the recent vector small-gain theorem given in [3]. Our results will extend the results obtained in [5] in several ways:

- we will show that approximate predictor schemes can be utilized under appropriate assumptions for the non-delayed system (1.1),
- we will propose implementation schemes for the (approximate or not) predictor-based feedback and
- we will propose explicit approximate predictor schemes for globally Lipschitz systems of the form (1.1), which are not necessarily feedforward systems.

Particularly, we will show that our main results (Theorem 2.2, Theorem 2.3, Corollary 3.4 and Corollary 3.6) can be applied to nonlinear triangular systems of the form:

$$\dot{x}_i(t) = f_i(x_1(t),...,x_i(t)) + x_{i+1}(t) \quad , \quad i = 1,...,n-1$$
$$\dot{x}_n(t) = f_n(x(t)) + u(t-r) \quad (1.3)$$
$$x(t) = (x_1(t),...,x_n(t))' \in \Re^n \; , \; u(t) \in \Re$$

where $f_i \in C^2(\Re^i; \Re)$ ($i = 1,...,n$) are globally Lipschitz functions.

The structure of the paper is as follows: Section 2 contains results which show that stabilization of (1.1) can be achieved by means of approximate predictors. Section 3 is devoted to the presentation of a systematic construction methodology of approximate predictors for globally Lipschitz systems. In Section 4 a simple example is provided, which illustrates the use of approximate predictor schemes. Finally, in Section 5 we present the concluding remarks of the present work.

**Notations** Throughout this paper we adopt the following notations:
* For a vector $x \in \Re^n$ we denote by $|x|$ its usual Euclidean norm, by $x'$ its transpose and by $|A| := \sup\{|Ax| \, ; \, x \in \Re^n, |x| = 1\}$ the induced norm of a matrix $A \in \Re^{m \times n}$.
* $\Re^+$ denotes the set of non-negative real numbers.
* For the definition of the class of functions $KL$, see [4].
* By $C^j(A)$ ($C^j(A; \Omega)$), where $j \geq 0$ is a non-negative integer, we denote the class of functions (taking values in $\Omega$) that have continuous derivatives of order $j$ on $A$.
* Let $x:[a-r,b] \to \Re^n$ with $b > a \geq 0$ and $r \geq 0$. By $T_r(t)x$ we denote the "history" of $x$ from $t-r$ to $t$, i.e., $T_r(t)x := \{x(t+\theta) \, ; \, \theta \in [-r,0]\}$, for $t \in [a,b)$.
* Let $I \subseteq \Re^+ := [0,+\infty)$ be an interval. By $L^\infty(I;U)$ ($L^\infty_{loc}(I;U)$) we denote the space of measurable and (locally) essentially bounded functions $u(\cdot)$ defined on $I$ and taking values in $U \subseteq \Re^m$. For $x \in L^\infty([-r,0]; \Re^n)$ we define $\|x\|_r := \sup_{\theta \in [-r,0]} |x(\theta)|$. We will also use the notation $M_U$ for the space of measurable and locally essentially bounded functions $u: \Re^+ \to U$.
* A continuous mapping $f: C^0([-r,0]; \Re^n) \times U \to \Re^k$, where $U \subseteq \Re^m$, is said to be completely Lipschitz with respect to $(x,u) \in C^0([-r,0]; \Re^n) \times U$ if for every bounded set $S \subset C^0([-r,0]; \Re^n) \times U$ there exists $L \geq 0$ such that $|f(x,u) - f(y,v)| \leq L \|x-y\|_r + L|u-v|$ for all $(x,u) \in S$, $(y,v) \in S$.
* Let $U \subseteq \Re^m$ be a closed non-empty convex set. For every $w \in \Re^m$, $\Pr_U(w)$ denotes the projection of $w$ on $U$.

We will always assume that the mapping $f: \Re^n \times U \to \Re^n$ is locally Lipschitz. Consequently, for every $(x_0, u) \in \Re^n \times M_U$ system (1.2) admits a unique local solution with initial condition $x(0) = x_0 \in \Re^n$ and corresponding to input $u \in M_U$.

## 2. Stabilization by Means of Approximate Predictors

We start by presenting the assumptions for system (1.2). We say that a system of the form (1.2) is forward complete if for every $x_0 \in \Re^n$, $u \in M_U$ the solution $x(t)$ of (1.2) with initial condition $x(0) = x_0 \in \Re^n$ corresponding to input $u:[-r,0] \to U$ exists for all $t \geq 0$. Using the semigroup property it is clear that (1.2) is forward complete if and only if there exists $r > 0$ such that for every $x_0 \in \Re^n$, $u \in M_U$ the solution $x(t)$ of (1.2) with initial condition $x(0) = x_0 \in \Re^n$ corresponding to input $u:[-r,0] \to U$ exists for all $t \in [0,r]$. Our first assumption concerning system (1.2) is the following:

**(H1)** System (1.2) is forward complete.

Assumption (H1) is also a necessary condition for the global stabilization of system (1.1) with $r > 0$: indeed, the solution of (1.1) for $t \in [0,r]$ must exist for every initial condition $x(0) = x_0 \in \Re^n$ and arbitrary input $u:[-r,0] \to U$. Therefore it follows that (1.2) must necessarily be a forward complete system.

Let $\phi(t, x_0; u)$ denote the solution of (1.2) at time $t \geq 0$ with initial condition $x(0) = x_0 \in \Re^n$ corresponding to input $u \in M_U$. The reader should notice that the solution $x(t)$ of (1.1) with initial condition $x(0) = x_0 \in \Re^n$ and corresponding to input $u \in C^0([-r,+\infty); U)$ satisfies:

$$x(t) = \phi(t, x_0; \delta_{-r} u) \text{ and } x(t+\tau) = \phi(t, x(\tau); \delta_{\tau-r} u) \text{ for all } t, \tau \geq 0 \qquad (2.1)$$

where $\delta_\theta : L_{loc}^\infty([-r,+\infty); U) \to M_U$ is the shift operator defined by

$$(\delta_\theta u)(t) := u(t+\theta), \text{ for } t \geq 0 \qquad (2.2)$$

We will assume next that (1.2) is stabilizable.

**(H2)** There exists $k \in C^1(\Re^n; U)$ with $k(0) = 0$ such that $0 \in \Re^n$ is Globally Asymptotically Stable for system (1.2) with $u = k(x)$, i.e., there exists a function $\sigma \in KL$ such that for every $x_0 \in \Re^n$ the solution $x(t)$ of (1.2) with $u = k(x)$ and initial condition $x(0) = x_0 \in \Re^n$ satisfies the following inequality:

$$|x(t)| \leq \sigma(|x_0|, t), \ \forall t \geq 0 \qquad (2.3)$$

**Theorem 2.1 (see Krstic, [5]):** *Consider system (1.1) under hypotheses (H1), (H2). Then the feedback law $u(t) = k(\phi(r, x(t); \delta_{t-r} u))$ globally asymptotically stabilizes system (1.1), i.e., there exists $\tilde\sigma \in KL$ such that for every $(x_0, u_0) \in \Re^n \times L^\infty([-r,0]; U)$ the solution $(x(t), u(t))$ of (1.2) with $u(t) = k(\phi(r, x(t); \delta_{t-r} u))$ and initial condition $x(0) = x_0 \in \Re^n$, $T_r(0)u = u_0$ satisfies the following inequality:*

$$|x(t)| + |u(t)| \leq \sigma(|x_0| + \|u_0\|_r, t), \ \forall t \geq 0 \qquad (2.4)$$

**Remark on Theorem 2.1:** It is clear that the implementation of the feedback law $u(t) = k(\phi(r, x(t); \delta_{t-r} u))$ involves the solution of an integral equation. The integral equation $u(t) = k(\phi(r, x(t); \delta_{t-r} u))$ may be transformed to a differential equation under certain regularity assumptions or can be given implicitly by the solution of a system of first order hyperbolic partial differential equations (see [5]). For the practical implementation of the feedback law the knowledge of the mapping $x \to \phi(r, x; u)$ is crucial.

We will next address the problem of the knowledge of the mapping $x \to \phi(r, x; u)$ and the implementation of the feedback law $u(t) = k(\phi(r, x(t); \delta_{t-r}u))$ by means of an approximate predictor scheme. Our hypotheses concerning systems (1.1) and (1.2) follow.

**(A1)** $f: \Re^n \times U \to \Re^n$ is locally Lipschitz with respect to $x \in \Re^n$ and there exists a constant $L \geq 0$ such that

$$x'f(x,u) \leq L|x|^2 + L|u|^2, \quad \forall x \in \Re^n, \forall u \in U$$

Hypothesis (A1) is a growth condition which guarantees hypothesis (H1). Indeed, by utilizing the function $V(x) = \frac{1}{2}|x|^2$, hypothesis (A1) implies that the derivative of $V(\phi(t, x_0; u))$ satisfies $\frac{d}{dt}V(\phi(t, x_0; u)) \leq 2LV(\phi(t, x_0; u)) + L|u(t)|^2$ for every $(x_0, u) \in \Re^n \times M_U$ and for all $t \geq 0$ for which $\phi(t, x_0; u)$ exists. Direct integration of the previous differential inequality implies

$$|\phi(t, x_0; u)| \leq \exp(Lt)\left(|x_0| + \sup_{0 \leq s \leq t}|u(s)|\right) \tag{2.5}$$

A standard contradiction argument in conjunction with (2.5) guarantees that (2.5) holds for every $(x_0, u) \in \Re^n \times M_U$ and for all $t \geq 0$.

**(A2)** There exists $k \in C^1(\Re^n; U)$ being locally Lipschitz with $k(0) = 0$ such that system (1.2) with $u = \Pr_U(k(x) + v)$ is Input-to-State Stable from the input $v \in \Re^m$ with linear gain function, i.e., there exist a function $\sigma \in KL$ and a constant $\gamma \geq 0$ such that for every $x_0 \in \Re^n$, $v \in M_{\Re^m}$, $t \geq 0$ the solution $x(t)$ of (1.2) with $u = \Pr_U(k(x) + v)$, initial condition $x(0) = x_0$ corresponding to input $v \in M_{\Re^m}$ satisfies the following inequality for all $t \geq 0$:

$$|x(t)| \leq \max\left\{\sigma(|x_0|, t), \gamma \sup_{0 \leq \tau \leq t}|v(\tau)|\right\} \tag{2.6}$$

Moreover, there exists a constant $R \geq 0$ such that

$$|k(x)| \leq R|x|, \quad \forall x \in \Re^n \tag{2.7}$$

Hypothesis (A2) is a more demanding hypothesis than (H2). The notion of Input-to-State Stability used here is the notion introduced by Sontag in [10].

Finally, we proceed with our last assumption concerning systems (1.1) and (1.2).

**(A3)** There exist constants $a_1, a_2 \geq 0$, $G \geq 0$ and completely Lipschitz mappings $p: \Re^n \times C^0([-r, 0]; U) \to U$, $g: \Re^n \times C^0([-r, 0]; U) \to \Re^m$ satisfying the following inequalities for all $(x, u) \in \Re^n \times C^0([-r, 0]; U)$:

$$|k(\phi(r, x; \delta_{-r}u)) - p(x, u)| \leq \max\{a_1|x|, a_2\|u\|_r\} \tag{2.8}$$

$$|g(x, u)| + |p(x, u)| \leq G|x| + G\|u\|_r \tag{2.9}$$

Moreover, for every $(x, u) \in \Re^n \times C^0([-r, +\infty); U)$ the solution $x(t)$ of (1.2) with initial condition $x(0) = x_0$ corresponding to $u \in C^0([-r, +\infty); U)$ satisfies $\frac{d}{dt}p(x(t), T_r(t)u) = g(x(t), T_r(t)u)$ for all $t \geq 0$ for which the solution exists.

Hypothesis (A3) introduces the mapping $p:\Re^n \times C^0([-r,0];U) \to U$, which approximates the stabilizing mapping $k(\phi(r,x;\delta_{-r}u))$. Indeed, the reader should notice that by virtue of (2.1) and (2.8) the solution $x(t)$ of (1.2) with initial condition $x(0) = x_0$ corresponding to $u \in C^0([-r,+\infty);U)$ satisfies:

$$|k(\phi(r,x(t);\delta_{t-r}u)) - p(x(t),T_r(t)u)| \leq \max\{a_1|x(t)|, a_2\|T_r(t)u\|_r\}, \quad \forall t \geq 0 \tag{2.10}$$

$$|k(x(t)) - p(x(t-r),T_r(t-r)u)| \leq \max\{a_1|x(t-r)|, a_2\|T_r(t-r)u\|_r\}, \quad \forall t \geq r \tag{2.11}$$

$$\frac{d}{dt}p(x(t),T_r(t)u) = g(x(t),T_r(t)u), \quad \forall t \geq 0 \tag{2.12}$$

Because of inequalities (2.10), (2.11), we will call the mapping $p:\Re^n \times C^0([-r,0];U) \to U$ an "approximate predictor" for (1.2). The constants $a_1, a_2 \geq 0$ determine how well the approximate predictor approximates the exact predictor scheme $k(\phi(r,x;\delta_{-r}u))$. For $a_1 = a_2 = 0$, we obtain the exact predictor, i.e., it holds that $k(\phi(r,x;\delta_{-r}u)) = p(x,u)$. Notice that for the exact predictor scheme, identity (2.12) holds with $g(x,u) = \nabla k(\phi(r,x;\delta_{-r}u))f(\phi(r,x;\delta_{-r}u),u(0))$.

We are now ready to state our first main result.

**Theorem 2.2:** *Consider systems (1.1) and (1.2) under hypotheses (A1-3) and further assume that*

$$\gamma a_1 < 1, \quad a_2(1+\gamma R) < 1 \tag{2.13}$$

*Then for every $\mu > 0$, there exists $\tilde{\sigma} \in KL$ such that for every $(x_0,w_0) \in \Re^n \times C^0([-r,0];\Re^m)$ the solution $(x(t),w(t))$ of (1.1) with*

$$u(t) = \Pr_U(w(t)) \tag{2.14}$$

$$\dot{w}(t) = g(x(t),T_r(t)u) - \mu(w(t) - p(x(t),T_r(t)u)) \tag{2.15}$$

*with initial condition $x(0) = x_0$, $T_r(0)w = w_0$ satisfies the estimate:*

$$|x(t)| + |w(t)| \leq \tilde{\sigma}(|x_0| + \|w_0\|_r, t), \quad \forall t \geq 0 \tag{2.16}$$

*i.e., the dynamic feedback law (2.14), (2.15) achieves global stabilization of system (1.1). Moreover, if there exist constants $M, \omega > 0$ such that the following estimate holds instead of (2.6)*

$$|x(t)| \leq \max\left\{M\exp(-\omega t)|x_0|, \gamma \sup_{0\leq\tau\leq t}\exp(-\omega(t-\tau))|v(\tau)|\right\}, \quad \forall t \geq 0 \tag{2.6'}$$

*then for every $\mu > 0$ there exist $\tilde{M}, \tilde{\omega} > 0$ such that for every $(x_0,w_0) \in \Re^n \times C^0([-r,0];\Re^m)$ the solution $(x(t),w(t))$ of (1.1), (2.14), (2.15) with initial condition $x(0) = x_0$, $T_r(0)w = w_0$ satisfies estimate (2.16) with $\tilde{\sigma}(s,t) := \tilde{M}\exp(-\tilde{\omega}t)s$, i.e., the dynamic feedback law (2.14), (2.15) achieves global exponential stabilization of system (1.1).*

**Remarks on Theorem 2.2:** Theorem 2.2 shows that approximate predictors can be used for the stabilization of system (1.1) provided that the approximation is sufficiently close to the exact predictor scheme. Moreover, Theorem 2.2 shows that the stabilizing feedback can be implemented as a dynamic feedback law; there is no need to solve integral equations. In general the stabilizing feedback law will involve distributed delays. To see this, notice that the classical Smith predictor (see [9]) for the linear system $f(x,u) = Ax + Bu$ with $k(x) = k'x$, $U = \Re^m$, where $A \in \Re^{n \times n}$, $B \in \Re^{n \times m}$, $k \in \Re^{n \times m}$ are constant matrices and $A + Bk'$ is Hurwitz, satisfies hypotheses (A1), (A2), (A3) with $a_1 = a_2 = 0$ and

$$\phi(r, x; u) := \exp(Ar)x + \int_0^r \exp(A(r-s))Bu(s)ds$$

$$p(x, u) := k'\exp(Ar)x + \int_0^r k'\exp(A(r-s))Bu(s-r)ds$$

$$g(x, u) := k'\exp(Ar)Ax + \int_0^r k'A\exp(A(r-s))Bu(s-r)d\theta + k'Bu(0)$$

Consequently, Theorem 2.2 guarantees that the closed-loop system with the dynamic distributed-delay feedback (2.14), (2.15), i.e., the system

$$\dot{x}(t) = Ax(t) + Bu(t-r)$$

$$\dot{u}(t) = k'\exp(Ar)(A + \mu I_n)Ax(t) + \int_0^r k'(A + \mu I_n)\exp(A(r-s))Bu(t+s-r)ds + (k'B - \mu I_m)u(t)$$

where $I_n \in \Re^{n \times n}$ denotes the identity matrix, is exponentially stable for all $\mu > 0$. Another important case where Theorem 2.2 is directly applicable is the case where hypothesis (A2) holds for certain $k \in C^2(\Re^n; U)$ satisfying

$$|\nabla k(x)| \leq R, \ \forall x \in \Re^n \tag{2.17}$$

and there exist constants $L_1, L_2 > 0$ such that the locally Lipschitz mapping $f(x, u)$ satisfies the growth condition

$$|f(x, u)| \leq L_1|x| + L_2|u|, \ \forall x \in \Re^n, u \in U \tag{2.18}$$

Indeed, in this case hypothesis (A1) automatically holds. We can distinguish two important cases where Theorem 2.2 is directly applicable:

- the case of the exact predictor scheme $p(x, u) = k(\phi(r, x; \delta_{-r}u))$, $g(x, u) = \nabla k(\phi(r, x; \delta_{-r}u))f(\phi(r, x; \delta_{-r}u), u(0))$. Indeed, utilizing (2.17), (2.18) and (2.5) with $L = L_1 + L_2$, it can be shown that hypothesis (A3) holds in this case with $a_1 = a_2 = 0$ and appropriate $G > 0$. Theorem 2.2 implies that the dynamic feedback law

$$u(t) = \Pr_U(w(t)) \tag{2.19}$$

$$\dot{w}(t) = \nabla k(\phi(r, x(t); \delta_{-r}T_r(t)u))f(\phi(r, x(t); \delta_{-r}T_r(t)u), u(t)) - \mu(w(t) - k(\phi(r, x(t); \delta_{-r}T_r(t)u))) \tag{2.20}$$

achieves global stabilization of system (1.1).

- the no predictor case, i.e., the case where $p(x, u) = k(x)$ and $g(x, u) = \nabla k(x)f(x, u(-r))$. The reader should notice that in this case there is no prediction and the stabilizing feedback law for (1.2) is used without any modification. Utilizing (2.17) and (2.18) we obtain for all $\varepsilon, \lambda > 0$, $(x, u) \in \Re^n \times C^0([-r, 0]; U)$

$$\sup_{0 \leq s \leq r} |\phi(s, x; \delta_{-r}u)| \leq \exp\left(\frac{2+\varepsilon}{2}L_1 r\right)|x| + \frac{L_2}{\varepsilon L_1}\sqrt{\frac{\varepsilon}{2+\varepsilon}}\sqrt{\exp((2+\varepsilon)L_1 r) - 1}\|u\|_r$$

and

$$|k(\phi(r,x;\delta_{-r}u)) - p(x,u)| \leq R|\phi(r,x;\delta_{-r}u) - x| \leq R\left|\int_0^r f(\phi(s,x;\delta_{-r}u), u(s-r))ds\right|$$

$$\leq R\int_0^r |f(\phi(s,x;\delta_{-r}u), u(s-r))|ds \leq R\int_0^r \left(L_1|\phi(s,x;\delta_{-r}u)| + L_2|u(s-r)|\right)ds$$

$$\leq RL_1 r \exp\left(\frac{2+\varepsilon}{2}L_1 r\right)|x| + R L_2 r\left(\frac{1}{\varepsilon}\sqrt{\frac{\varepsilon}{2+\varepsilon}}\sqrt{\exp((2+\varepsilon)L_1 r)-1} + 1\right)\|u\|_r$$

$$\leq R r \max\left\{L_1(1+\lambda)\exp\left(\frac{2+\varepsilon}{2}L_1 r\right)|x|, L_2\left(1+\frac{1}{\lambda}\right)\left(\frac{1}{\varepsilon}\sqrt{\frac{\varepsilon}{2+\varepsilon}}\sqrt{\exp((2+\varepsilon)L_1 r)-1} + 1\right)\|u\|_r\right\}$$

It follows that hypothesis (A3) holds. By virtue of Theorem 2.2, the dynamic feedback law

$$u(t) = \Pr_U(w(t)) \tag{2.21}$$

$$\dot{w}(t) = \nabla k(x(t)) f(x(t), u(t-r)) - \mu(w(t) - k(x(t))) \tag{2.22}$$

achieves global stabilization of system (1.1), provided that the delay $r > 0$ is sufficiently small. More specifically, the above inequalities show that global stabilization of system (1.1) is achieved provided that there exists $\varepsilon > 0$ such that

$$(1+\gamma R)R\, rL_2\left(\frac{1}{\varepsilon}\sqrt{\frac{\varepsilon}{2+\varepsilon}}\sqrt{\exp((2+\varepsilon)L_1 r)-1}+1\right) + R\, rL_1\gamma \exp\left(\frac{2+\varepsilon}{2}L_1 r\right) < 1$$

where $\gamma \geq 0$ is the constant involved in hypothesis (A2).

**Proof of Theorem 2.2:** Since $U \subseteq \Re^m$ is closed and convex with $0 \in U$ the following inequalities will be used repeatedly in the proof:

$$|\Pr_U(w)| \leq |w| \text{ and } |\Pr_U(w) - \Pr_U(v)| \leq |w-v|, \text{ for all } w,v \in \Re^m$$

Exploiting hypothesis (A1) and the linear growth condition (2.9), it can be shown (using the functional $V(x,w) = |x|^2 + \sup_{-r \leq s \leq 0} |w(s)|^2$) that the solution of (1.1), (2.14), (2.15) with initial condition $x(0) = x_0$, $T_r(0)w = w_0$ exists for all $t \geq 0$ and satisfies the inequality:

$$|x(t)| + \|T_r(t)w\|_r \leq B\exp(\sigma t)\left(|x_0| + \|w_0\|_r\right), \quad \forall t \geq 0 \tag{2.23}$$

for certain constants $B, \sigma > 0$. Differential equation (2.15) and hypothesis (A2) imply that the following inequalities hold for the solution of (1.1), (2.14), (2.15) with initial condition $x(0) = x_0$, $T_r(0)w = w_0$:

$$|w(t) - p(x(t), T_r(t)u)| \leq \exp(-\mu t)|w(0) - p(x_0, T_r(0)u)|, \quad \forall t \geq 0 \tag{2.24}$$

$$|x(t)| \leq \max\left\{\sigma(|x_0|, t), \gamma \sup_{0 \leq s \leq t}|w(s-r) - k(x(s))|\right\}, \quad \forall t \geq 0 \tag{2.25}$$

Using (2.5) and (2.7) we obtain for $t \in [0, r)$:

$$|w(t-r) - k(x(t))| \leq |w(t-r)| + R|x(t)| \leq \|w_0\|_r + R|x(t)|$$
$$\leq \|w_0\|_r + R\exp(Lr)\left(|x_0| + \sup_{-r\leq s\leq 0}|u(s)|\right)$$
$$\leq (1 + R\exp(Lr))\|w_0\|_r + R\exp(Lr)|x_0| \qquad (2.26)$$
$$\leq (1 + R\exp(Lr))\exp(-\mu(t-r))(|x_0| + \|w_0\|_r)$$

Using (2.9), (2.11) and (2.24) we obtain for $t \geq r$:

$$|w(t-r) - k(x(t))| \leq$$
$$\leq |w(t-r) - p(x(t-r), T_r(t-r)u)| + |p(x(t-r), T_r(t-r)u) - k(x(t))|$$
$$\leq \exp(-\mu(t-r))|w(0) - p(x_0, T_r(0)u)| + \max\{a_1|x(t-r)|, a_2\|T_r(t-r)u\|_r\} \qquad (2.27)$$
$$\leq (1+G)\exp(-\mu(t-r))(\|w_0\|_r + |x_0|) + \max\{a_1|x(t-r)|, a_2\|T_r(t-r)u\|_r\}$$
$$\leq (1+G)\exp(-\mu(t-r))(\|w_0\|_r + |x_0|) + \max\{a_1|x(t-r)|, a_2\|T_r(t-r)w\|_r\}$$

Combining (2.26) and (2.27) we conclude that there exists a constant $Q > 0$ such that the following inequality holds for all $t \geq 0$ and $\varepsilon > 0$:

$$|w(t-r) - k(x(t))| \leq \max\left\{Q(1+\varepsilon^{-1})\exp(-\mu t)(\|w_0\|_r + |x_0|), a_1(1+\varepsilon)\sup_{0\leq s\leq t}|x(s)|, a_2(1+\varepsilon)\sup_{0\leq s\leq t}|w(s-r)|\right\} \qquad (2.28)$$

On the other hand, using (2.7) we conclude that the following inequalities hold for all $t \geq 0$ and $\lambda > 0$:

$$|w(t-r)| \leq |w(t-r) - k(x(t))| + |k(x(t))|$$
$$\leq |w(t-r) - k(x(t))| + R|x(t)|$$
$$\leq \sup_{0\leq s\leq t}|w(s-r) - k(x(s))| + R\sup_{0\leq s\leq t}|x(s)| \qquad (2.29)$$
$$\leq \max\left\{(1+\lambda)\sup_{0\leq s\leq t}|w(s-r) - k(x(s))|, (1+\lambda^{-1})R\sup_{0\leq s\leq t}|x(s)|\right\}$$

Using (2.23), (2.25), (2.28) and (2.29) in conjunction with the vector small-gain theorem in [3], we conclude that there exists $\tilde{\sigma} \in KL$ such that (2.16) holds, provided that there exist $\varepsilon, \lambda > 0$ satisfying the following inequalities:

$$(1+\varepsilon)\gamma a_1 < 1, \quad (1+\varepsilon)a_2(1+\lambda^{-1})\gamma R < 1, \quad a_2(1+\varepsilon)(1+\lambda) < 1$$

The reader should notice that inequalities (2.13) guarantee the existence $\varepsilon, \lambda > 0$ such that the above inequalities hold.

To finish the proof, consider the case where (2.6') holds for certain constants $M, \omega > 0$. Let $0 < \tilde{\omega} \leq \min\{\omega, \mu\}$ sufficiently small such that

$$(1+\varepsilon)\gamma a_1 \exp(\tilde{\omega} r) < 1, \quad (1+\varepsilon)a_2\exp(\tilde{\omega} r)(1+\lambda^{-1})\gamma R < 1, \quad a_2\exp(\tilde{\omega} r)(1+\varepsilon)(1+\lambda) < 1 \qquad (2.30)$$

for certain $\varepsilon, \lambda > 0$. Again the existence of appropriate $\varepsilon, \lambda > 0$ and sufficiently small $\tilde{\omega} > 0$ is guaranteed by (2.13). Inequality (2.6') gives:

$$\exp(\tilde{\omega}t)|x(t)| \leq \max\left\{M|x_0|, \gamma \sup_{0\leq s\leq t}\exp(\tilde{\omega} s)|w(s-r) - k(x(s))|\right\}, \quad \forall t \geq 0 \qquad (2.31)$$

Using (2.9), (2.11) and (2.24) we obtain for $t \geq r$ and $\tilde{\omega} \leq \mu$:

$$\exp(\widetilde{\omega}t)|w(t-r)-k(x(t))| \leq$$
$$\leq \exp(\widetilde{\omega}t)|w(t-r)-p(x(t-r),T_r(t-r)u)| + \exp(\widetilde{\omega}t)|p(x(t-r),T_r(t-r)u)-k(x(t))|$$
$$\leq \exp(\widetilde{\omega}t)\exp(-\mu(t-r))|w(0)-p(x_0,T_r(0)u)| + \max\{a_1\exp(\widetilde{\omega}t)|x(t-r)|, a_2\exp(\widetilde{\omega}t)\|T_r(t-r)u\|_r\}$$
$$\leq (1+G)\exp(\mu r)(\|w_0\|_r + |x_0|) + \max\{a_1\exp(\widetilde{\omega}t)|x(t-r)|, a_2\exp(\widetilde{\omega}t)\|T_r(t-r)u\|_r\}$$
$$\leq (1+G)\exp(\mu r)(\|w_0\|_r + |x_0|) + \max\left\{a_1\exp(\widetilde{\omega}r)\exp(\widetilde{\omega}(t-r))|x(t-r)|, a_2\exp(\widetilde{\omega}t)\sup_{-r\leq s\leq 0}|w(t-r+s)|\right\}$$
$$\leq (1+G)\exp(\mu r)(\|w_0\|_r + |x_0|) + \max\left\{a_1\exp(\widetilde{\omega}r)\exp(\widetilde{\omega}(t-r))|x(t-r)|, a_2\exp(2\widetilde{\omega}r)\sup_{-r\leq s\leq 0}\exp(\widetilde{\omega}(t-r+s))|w(t-r+s)|\right\}$$

Combining (2.26) and the above inequality we conclude that there exists a constant $\widetilde{Q} > 0$ such that the following inequality holds for all $t \geq 0$:

$$\exp(\widetilde{\omega}t)|w(t-r)-k(x(t))| \leq \max\left\{\begin{array}{l}\widetilde{Q}(1+\varepsilon^{-1})(\|w_0\|_r + |x_0|), a_1\exp(\widetilde{\omega}r)(1+\varepsilon)\sup_{0\leq s\leq t}\exp(\widetilde{\omega}s)|x(s)|, \\ a_2\exp(2\widetilde{\omega}r)(1+\varepsilon)\sup_{0\leq s\leq t}\exp(\widetilde{\omega}(s-r))|w(s-r)|\end{array}\right\} \quad (2.32)$$

On the other hand, using (2.7) we conclude that the following inequalities hold for all $t \geq 0$:

$$\exp(\widetilde{\omega}(t-r))|w(t-r)| \leq \exp(\widetilde{\omega}(t-r))|w(t-r)-k(x(t))| + \exp(\widetilde{\omega}(t-r))|k(x(t))|$$
$$\leq |w(t-r)-k(x(t))|\exp(\widetilde{\omega}(t-r)) + R|x(t)|\exp(\widetilde{\omega}(t-r))$$
$$\leq \exp(-\widetilde{\omega}r)\sup_{0\leq s\leq t}|w(s-r)-k(x(s))|\exp(\widetilde{\omega}s) + R\exp(-\widetilde{\omega}r)\sup_{0\leq s\leq t}|x(s)|\exp(\widetilde{\omega}s) \quad (2.33)$$
$$\leq \max\left\{(1+\lambda)\exp(-\widetilde{\omega}r)\sup_{0\leq s\leq t}|w(s-r)-k(x(s))|\exp(\widetilde{\omega}s), (1+\lambda^{-1})R\exp(-\widetilde{\omega}r)\sup_{0\leq s\leq t}|x(s)|\exp(\widetilde{\omega}s)\right\}$$

Combining (2.32) and (2.33) we obtain for all $t \geq 0$:

$$\sup_{0\leq s\leq t}\exp(\widetilde{\omega}s)|w(s-r)-k(x(s))| \leq \max\left\{\begin{array}{l}\widetilde{Q}(1+\varepsilon^{-1})(\|w_0\|_r + |x_0|), a_1\exp(\widetilde{\omega}r)(1+\varepsilon)\sup_{0\leq s\leq t}\exp(\widetilde{\omega}s)|x(s)|, \\ a_2\exp(\widetilde{\omega}r)(1+\varepsilon)(1+\lambda)\sup_{0\leq s\leq t}|w(s-r)-k(x(s))|\exp(\widetilde{\omega}s), \\ a_2\exp(\widetilde{\omega}r)(1+\varepsilon)(1+\lambda^{-1})R\sup_{0\leq s\leq t}|x(s)|\exp(\widetilde{\omega}s)\end{array}\right\} \quad (2.34)$$

Since $a_2\exp(\widetilde{\omega}r)(1+\varepsilon)(1+\lambda) < 1$ (recall (2.30)), inequality (2.34) is simplified in the following way:

$$\sup_{0\leq s\leq t}\exp(\widetilde{\omega}s)|w(s-r)-k(x(s))| \leq \max\left\{\begin{array}{l}\widetilde{Q}(1+\varepsilon^{-1})(\|w_0\|_r + |x_0|), \\ \exp(\widetilde{\omega}r)(1+\varepsilon)\max\{a_1, a_2(1+\lambda^{-1})R\}\sup_{0\leq s\leq t}|x(s)|\exp(\widetilde{\omega}s)\end{array}\right\} \quad (2.35)$$

Inequality (2.35) in conjunction with inequality (2.31) gives for all $t \geq 0$:

$$\sup_{0\leq s\leq t}\exp(\widetilde{\omega}s)|w(s-r)-k(x(s))| \leq \max\left\{\begin{array}{l}\widetilde{Q}(1+\varepsilon^{-1})(\|w_0\|_r + |x_0|), \\ \exp(\widetilde{\omega}r)(1+\varepsilon)\max\{a_1, a_2(1+\lambda^{-1})R\}M|x_0|, \\ \exp(\widetilde{\omega}r)(1+\varepsilon)\max\{a_1, a_2(1+\lambda^{-1})R\}\gamma\sup_{0\leq s\leq t}\exp(\widetilde{\omega}s)|w(s-r)-k(x(s))|\end{array}\right\}$$

Since $\exp(\widetilde{\omega}r)(1+\varepsilon)\max\{a_1, a_2(1+\lambda^{-1})R\}\gamma < 1$ (recall (2.30)), we obtain:

$$\sup_{0\leq s\leq t} \exp(\tilde{\omega}s)|w(s-r)-k(x(s))| \leq \max\left\{\begin{array}{l} \tilde{Q}(1+\varepsilon^{-1})(\|w_0\|_r +|x_0|), \\ \exp(\tilde{\omega}r)(1+\varepsilon)\max\{a_1,a_2(1+\lambda^{-1})R\}M|x_0| \end{array}\right\} \quad (2.36)$$

Inequality (2.36) in conjunction with (2.31) gives:

$$\sup_{0\leq s\leq t} \exp(\tilde{\omega}s)|x(s)| \leq \max\left\{\begin{array}{l} M|x_0|, \gamma\tilde{Q}(1+\varepsilon^{-1})(\|w_0\|_r +|x_0|), \\ \gamma\exp(\tilde{\omega}r)(1+\varepsilon)\max\{a_1,a_2(1+\lambda^{-1})R\}M|x_0| \end{array}\right\}, \quad \forall t\geq 0 \quad (2.37)$$

Finally, from (2.36), (2.37) and (2.33) we get:

$$\exp(\tilde{\omega}(t-r))|w(t-r)| \leq P(\|w_0\|_r +|x_0|), \quad \forall t\geq 0 \quad (2.38)$$

for certain appropriate constant $P>0$. Inequalities (2.37) and (2.38) imply that there exist $\tilde{M},\tilde{\omega}>0$ such that for every $(x_0,w_0)\in\Re^n\times C^0([-r,0];\Re^m)$ the solution $(x(t),w(t))$ of (1.1), (2.14), (2.15) with initial condition $x(0)=x_0$, $T_r(0)w=w_0$ satisfies estimate (2.16) with $\tilde{\sigma}(s,t):=\tilde{M}\exp(-\tilde{\omega}t)s$.

The proof is complete. ◁

We finish this section by providing an additional result on approximate predictors. Since the formulae for the mappings $p:\Re^n\times C^0([-r,0];U)\to U$ and $g:\Re^n\times C^0([-r,0];U)\to\Re^m$ involved in hypothesis (A3) are usually complicated (see next section), the following result helps for the simplification of the formulae at the cost of an additional approximation.

**Theorem 2.3:** *Consider systems (1.1) and (1.2) under hypotheses (A1-2) and further assume that the following hypothesis holds:*

**(A4)** There exist constants $a_1,a_2\geq 0$, $G\geq 0$ and completely Lipschitz mappings $p:\Re^n\times C^0([-r,0];U)\to U$, $g:\Re^n\times C^0([-r,0];U)\to\Re^m$ satisfying the following inequalities for all $(x,u)\in\Re^n\times C^0([-r,0];U)$:

$$\max\{|k(\phi(r,x;\delta_{-r}u))-p(x,u)|,|\nabla k(\phi(r,x;\delta_{-r}u))f(\phi(r,x;\delta_{-r}u),u(0))-g(x,u)|\}\leq \max\{a_1|x|,a_2\|u\|_r\} \quad (2.39)$$

$$|g(x,u)|+|p(x,u)|\leq G|x|+G\|u\|_r \quad (2.40)$$

$$\gamma a_1 <1, \quad a_2(1+\gamma R)<1 \quad (2.41)$$

*Then for every $\mu>0$ satisfying*

$$\gamma a_1\left(1+\frac{1}{\mu}\right)<1, \quad a_2\left(1+\frac{1}{\mu}\right)(1+\gamma R)<1 \quad (2.42)$$

*there exists $\tilde{\sigma}\in KL$ such that for every $(x_0,w_0)\in\Re^n\times C^0([-r,0];\Re^m)$ the solution $(x(t),w(t))$ of (1.1) with (2.14), (2.15) and initial condition $x(0)=x_0$, $T_r(0)w=w_0$ satisfies estimate (2.16), i.e., the dynamic feedback law (2.14), (2.15) achieves global stabilization of system (1.1). Moreover, if there exist constants $M,\omega>0$ such that (2.6') holds instead of (2.6), then for every $\mu>0$ satisfying (2.42) there exist $\tilde{M},\tilde{\omega}>0$ such that for every $(x_0,w_0)\in\Re^n\times C^0([-r,0];\Re^m)$ the solution $(x(t),w(t))$ of (1.1), (2.14), (2.15) with initial condition $x(0)=x_0$, $T_r(0)w=w_0$ satisfies estimate (2.16) with $\tilde{\sigma}(s,t):=\tilde{M}\exp(-\tilde{\omega}t)s$, i.e., the dynamic feedback law (2.14), (2.15) achieves global exponential stabilization of system (1.1).*

**Proof:** The proof is exactly the same with the proof of Theorem 2.2 except of the estimate for the quantity $|w(t-r) - k(x(t))|$. By virtue of inequality (2.39) and noticing that

$$\frac{d}{dt} k(\phi(r, x(t); \delta_{-r} T_r(t)u)) = \nabla k(\phi(r, x(t); \delta_{-r} T_r(t)u)) f(\phi(r, x(t); \delta_{-r} T_r(t)u), u(t))$$

we obtain

$$\left| \frac{d}{dt} k(\phi(r, x(t); \delta_{-r} T_r(t)u)) - g(x(t), T_r(t)u) \right| \leq \max\{a_1 |x(t)|, a_2 \|T_r(t)u\|_r\}$$

Integrating (2.15) and using (2.39) and the above inequality, we obtain for all $t \geq 0$:

$$\begin{aligned} &|w(t) - k(\phi(r, x(t); \delta_{-r} T_r(t)u))| \leq \\ &\exp(-\mu t)|w(0) - k(\phi(r, x(0); \delta_{-r} T_r(0)u))| + \left(1 + \frac{1}{\mu}\right) \max\left\{ a_1 \sup_{0 \leq \tau \leq t} |x(\tau)|, a_2 \sup_{0 \leq \tau \leq t} \|T_r(\tau)u\|_r \right\} \end{aligned} \quad (2.43)$$

Combining (2.26), (2.43), (2.5), (2.7) we conclude that there exists a constant $Q > 0$ such that the inequality (2.28) holds for all $t \geq 0$ and $\varepsilon > 0$ with $a_1, a_2$ replaced by $a_1(1 + \mu^{-1}), a_2(1 + \mu^{-1})$.

Similar changes are needed for the case of exponential stability. Details are left to the reader. ◁

## 3. Approximate Predictors for Globally Lipschitz Nonlinear Systems

In this section we show how we can construct approximate predictors for globally Lipschtz systems, i.e., systems for which there exists a constant $L \geq 0$ satisfying

$$|f(x, u) - f(y, u)| \leq L|x - y|, \quad \forall x, y \in \Re^n, \forall u \in U \quad (3.1a)$$

$$|f(x, u)| \leq L|x| + L|u|, \quad \forall x \in \Re^n, \forall u \in U \quad (3.1b)$$

Particularly, we will show that the solution map for system (1.2) under (3.1a,b) can be approximated by successive approximations.

Let $u \in L^\infty([0,T]; U)$ be arbitrary and define the operator $P_{T,u} : C^0([0,T]; \Re^n) \to C^0([0,T]; \Re^n)$ by

$$(P_{T,u} x)(t) = x(0) + \int_0^t f(x(\tau), u(\tau)) d\tau, \text{ for } t \in [0, T] \quad (3.2)$$

and we denote by $P_{T,u}^l = \underbrace{P_{T,u} \ldots P_{T,u}}_{l \text{ times}}$ for every integer $l \geq 1$. The following facts hold for the operator $P_{T,u} : C^0([0,T]; \Re^n) \to C^0([0,T]; \Re^n)$.

Fact I

$$\max_{0 \leq t \leq T} |(P_{T,u} x)(t) - (P_{T,u} y)(t)| \leq |x(0) - y(0)| + LT \max_{0 \leq \tau \leq T} |x(\tau) - y(\tau)|, \quad \forall x, y \in C^0([0,T]; \Re^n) \quad (3.3)$$

The above fact is a direct consequence of (3.1a).

Fact II

For all $x \in C^0([0,T]; \Re^n)$ the following implication holds:

$$\text{If } LT < 1 \text{ then } \max_{0 \leq t \leq T}\left|(P_{T,u}^l x)(t) - \phi(t, x(0); u))\right| \leq \frac{(LT)^l}{1-LT} \max_{0 \leq \tau \leq T}\left|x(0) + \int_0^t f(x(\tau), u(\tau))d\tau - x(t)\right| \quad (3.4)$$

The proof of the above fact follows closely the proof of Banach's fixed point theorem: first we show (by induction and using Fact I) that

$$\max_{0 \leq t \leq T}\left|(P_{T,u}^l x)(t) - (P_{T,u}^{l-1} x)(t)\right| \leq (LT)^{l-1} \max_{0 \leq \tau \leq T}\left|x(0) + \int_0^t f(x(\tau), u(\tau))d\tau - x(t)\right|$$

$$\text{for all } l \geq 1 \text{ and } x \in C^0([0,T]; \Re^n) \quad (3.5)$$

Then we proceed by estimating the quantity $\max_{0 \leq t \leq T}\left|(P_{T,u}^m x)(t) - x(t)\right|$ by using (3.5) and the inequality

$$\max_{0 \leq t \leq T}\left|(P_{T,u}^m x)(t) - x(t)\right| \leq \max_{0 \leq t \leq T}\left|(P_{T,u}^m x)(t) - (P_{T,u}^{m-1} x)(t)\right| + \ldots + \max_{0 \leq t \leq T}\left|(P_{T,u}^1 x)(t) - x(t)\right|$$

$$\leq \left[(LT)^{m-1} + \ldots + 1\right] \max_{0 \leq \tau \leq T}\left|x(0) + \int_0^t f(x(\tau), u(\tau))d\tau - x(t)\right| =$$

$$= \frac{1-(LT)^m}{1-LT} \max_{0 \leq \tau \leq T}\left|x(0) + \int_0^t f(x(\tau), u(\tau))d\tau - x(t)\right|$$

Replacing $x$ in the above inequality with $P_{T,u}^l x$ and using (3.5), we get for all $x \in C^0([0,T]; \Re^n)$:

$$\max_{0 \leq t \leq T}\left|(P_{T,u}^{m+l} x)(t) - (P_{T,u}^l x)(t)\right| \leq \frac{1-(LT)^m}{1-LT} \max_{0 \leq t \leq T}\left|(P_{T,u}^{l+1} x)(t) - (P_{T,u}^l x)(t)\right|$$

$$\leq \frac{1-(LT)^m}{1-LT}(LT)^l \max_{0 \leq \tau \leq T}\left|x(0) + \int_0^t f(x(\tau), u(\tau))d\tau - x(t)\right|$$

Finally, we notice that $\lim_{m \to \infty} \max_{0 \leq t \leq T}\left|(P_{T,u}^m x)(t) - \phi(t, x(0); u))\right| = 0$, by virtue of Banach's fixed point theorem. By letting $m \to +\infty$ in the above inequality we obtain (3.4).

We next define the operators $G_T : \Re^n \to C^0([0,T]; \Re^n)$, $C_T : C^0([0,T]; \Re^n) \to \Re^n$ and $Q_{T,u}^l : \Re^n \to \Re^n$ for $l \geq 1$ by

$$(G_T x_0)(t) = x_0, \text{ for } t \in [0,T] \text{ and } C_T x = x(T) \quad (3.6)$$

$$Q_{T,u}^l = C_T P_{T,u}^l G_T \quad (3.7)$$

The following fact holds for the mapping $Q_{T,u}^l : \Re^n \to \Re^n$.

Fact III

For every $x, y \in \Re^n$ and $u \in L^\infty([0,T]; U)$ the following implication holds:

$$\text{If } LT < 1 \text{ then } \left|Q_{T,u}^l x - \phi(T, y; u))\right| \leq \frac{(LT)^{l+1}}{1-LT}\left(|x| + \sup_{0 \leq \tau \leq T}|u(\tau)|\right) + \exp(LT)|x - y| \quad (3.8)$$

Proof of Fact III: Implication (3.4) and definitions (3.6), (3.7) imply that for every $x \in \Re^n$ the following implication holds:

$$\text{If } LT < 1 \text{ then } \left| Q_{T,u}^l x - \phi(T,x;u)) \right| \leq \frac{(LT)^l}{1-LT} \max_{0 \leq \tau \leq T} \left| \int_0^t f(x,u(\tau))d\tau \right| \quad (3.9)$$

On the other hand inequality (3.1a) gives for all $x, y \in \Re^n$:

$$\left| \phi(t,y;u)) - \phi(t,x;u)) \right| \leq |x-y| + L \int_0^t \left| \phi(s,y;u)) - \phi(s,x;u)) \right| ds, \quad \forall t \in [0,T] \quad (3.10)$$

Application of Gronwall's lemma (see [4]) to inequality (3.10) gives for all $x, y \in \Re^n$:

$$\left| \phi(T,y;u)) - \phi(T,x;u)) \right| \leq \exp(LT)|x-y| \quad (3.11)$$

Implication (3.8) is a direct consequence of implication (3.9), inequality (3.11) and the inequality $\max_{0 \leq \tau \leq T} \left| \int_0^t f(x,u(\tau))d\tau \right| \leq LT|x| + LT \sup_{0 \leq \tau \leq T} |u(\tau)|$ which is a consequence of hypothesis (3.1b).

The reader should notice at this point that inequality (3.8) guarantees that the mapping $Q_{T,u}^l : \Re^n \to \Re^n$ approximates the solution map $\phi(T,x;u))$ if $LT < 1$. Moreover, the approximation error can be tuned to be "small" by allowing $l \geq 1$ to take large values. Finally, the mapping $Q_{T,u}^l : \Re^n \to \Re^n$ is easily computed. The following example illustrates this point.

**Example 3.1:** Consider the case $f(x,u) = a(x) + b(u)$, where $a : \Re^n \to \Re^n$ is a globally Lipschitz vector field with Lipschitz constant $L \geq 0$ and $b : U \to \Re^n$ is a continuous mapping satisfying the linear growth condition $|b(u)| \leq L|u|$ for all $u \in U$. In this case inequalities (3.1a,b) hold. Applying definitions (3.2), (3.6), (3.7) we get for all $x \in \Re^n$ and $u \in L^\infty([0,T];U)$:

$$Q_{T,u}^1 x = x + Ta(x) + \int_0^T b(u(\tau))d\tau \quad (3.12)$$

$$Q_{T,u}^2 x = x + \int_0^T a\left( x + \tau a(x) + \int_0^\tau b(u(s))ds \right) d\tau + \int_0^T b(u(\tau))d\tau \quad (3.13)$$

Fact III guarantees that if $LT < 1$ then the following inequality will hold:

$$\left| Q_{T,u}^2 x - \phi(T,x;u)) \right| \leq \frac{(LT)^3}{1-LT} \left( |x| + \sup_{0 \leq \tau \leq T} |u(\tau)| \right) \quad (3.14)$$

The reader should notice that it is easy to generate mappings $Q_{T,u}^l : \Re^n \to \Re^n$ with $l > 2$. This example will be continued. ◁

We next define the mapping $P_{l,q}^u : \Re^n \to \Re^n$ for arbitrary $u \in L^\infty([0,r];U)$. Let $l, q \geq 1$ be integers and $T = \frac{r}{q}$. We define for all $x \in \Re^n$:

$$P_{l,q}^u x = Q_{T,u_q}^l \ldots Q_{T,u_1}^l x \tag{3.15}$$

where $u_i(s) = u(s+(i-1)T)$, $i=1,\ldots,q$ for $s \in [0,T]$. Notice that $u_i \in L^\infty([0,T];U)$ for $i=1,\ldots,q$. For the operator $P_{l,q}^u : \Re^n \to \Re^n$ we are in a position to prove the following proposition. Its proof is provided at the Appendix.

**Proposition 3.2:** *Let $l,q$ be positive integers with $LT < 1$, where $T = \dfrac{r}{q}$. Suppose that inequalities (3.1a,b) hold. Then there exists a constant $K := K(q) \geq 0$, independent of $l$, such that for every $u \in L^\infty([0,r];U)$ and $x \in \Re^n$ the following inequality holds:*

$$\left| P_{l,q}^u x - \phi(r,x;u) \right| \leq K \frac{(LT)^{l+1}}{1-LT} \left( |x| + \sup_{0 \leq \tau \leq r} |u(\tau)| \right) \tag{3.16}$$

The reader should again notice that inequality (3.16) guarantees that the mapping $P_{l,q}^u : \Re^n \to \Re^n$ approximates the solution map $\phi(r,x;u))$ if $Lr < q$. Again, the approximation error can be guaranteed to be "sufficiently small" by allowing $l \geq 1$ to take large values. Finally, the mapping $P_{l,q}^u : \Re^n \to \Re^n$ can be easily computed.

**Example 3.3:** Consider again the case $f(x,u) = a(x) + b(u)$, where $a : \Re^n \to \Re^n$ is a globally Lipschitz vector field with Lipschitz constant $L \geq 0$ and $b : U \to \Re^n$ is a continuous mapping satisfying the linear growth condition $|b(u)| \leq L|u|$ for all $u \in U$. Applying definitions (3.15), (3.12), (3.13) we get for all $x \in \Re^n$ and $u \in L^\infty([0,T];U)$:

$$P_{1,2}^u x = x + \frac{r}{2} a(x) + \int_0^r b(u(\tau))d\tau + \frac{r}{2} a\left( x + \frac{r}{2} a(x) + \int_0^{r/2} b(u(\tau))d\tau \right)$$

$$P_{2,2}^u x = x_1 + \int_0^{r/2} a\left( x_1 + \tau a(x_1) + \int_{r/2}^{\tau+r/2} b(u(s))ds \right) d\tau + \int_{r/2}^r b(u(\tau))d\tau$$

$$x_1 = x + \int_0^{r/2} a\left( x + \tau a(x) + \int_0^\tau b(u(s))ds \right) d\tau + \int_0^{r/2} b(u(\tau))d\tau$$

Proposition 3.2 guarantees that if $Lr < 2$ then there exists a constant $K \geq 0$ such that the following inequality holds for all $l \geq 1$:

$$\left| P_{l,2}^u x - \phi(r,x;u)) \right| \leq K \frac{(LT)^{l+1}}{1-LT} \left( |x| + \sup_{0 \leq \tau \leq r} |u(\tau)| \right) \tag{3.17}$$

Although formulae for $P_{l,q}^u : \Re^n \to \Re^n$ are complicated for large $l,q$, the values for $P_{l,q}^u x$ can be provided through a simple algorithm. ◁

Finally, let $k \in C^1(\Re^n;U)$ be a mapping with $k(0) = 0$ and for which there exists a constant $R \geq 0$ such that:

$$\left| \nabla k(x) \right| \leq R, \quad \forall x \in \Re^n \tag{3.18}$$

and consider next the mapping $p_{l,q} : \Re^n \times C^0([-r,0];U) \to U$ defined by:

$$p_{l,q}(x,u) := k\left( P_{l,q}^{\delta_{-r} u} x \right) \tag{3.19}$$

Notice that inequality (3.1b) guarantees that inequality (2.5) holds with $L$ replaced by $\frac{1+\sqrt{2}}{2}L$. Proposition 3.2 and definition (3.19) in conjunction with inequality (3.18) and inequality (2.5) (with $L$ replaced by $\frac{1+\sqrt{2}}{2}L$) guarantee the existence of a constant $K := K(q) \geq 0$, independent of $l$, such that for every $(x,u) \in \Re^n \times C^0([-r,0];U)$ and $\varepsilon > 0$ the following inequality holds:

$$\left| p_{l,q}(x,u) - k(\phi(r,x;\delta_{-r}u)) \right| \leq \max\left\{ (1+\varepsilon)RK\frac{(LT)^{l+1}}{1-LT}|x|, (1+\varepsilon^{-1})RK\frac{(LT)^{l+1}}{1-LT}\|u\|_r \right\} \tag{3.20}$$

$$\left| p_{l,q}(x,u) \right| \leq R\left[ \exp\left(\frac{1+\sqrt{2}}{2}Lr\right) + K\frac{(LT)^{l+1}}{1-LT} \right]\left(|x| + \|u\|_r\right) \tag{3.21}$$

where $l, q$ are positive integers and $T = \frac{r}{q}$ with $Lr < q$. Therefore, the mapping $p_{l,q} : \Re^n \times C^0([-r,0];U) \to U$ is a natural candidate to satisfy the requirements of hypothesis (A3) in Section 2, i.e., to be an approximate predictor. Indeed, Theorem 2.2 allows us to obtain the following corollary.

**Corollary 3.4:** *Consider systems (1.1) and (1.2) under hypothesis (A2) and further assume that inequalities (3.1a,b), (3.18) hold. Let $p_{l,q} : \Re^n \times C^0([-r,0];U) \to U$ be the mapping defined by (3.19) for positive integers $l, q$ with $Lr < q$. Moreover, assume that $p_{l,q} : \Re^n \times C^0([-r,0];U) \to U$ is completely Lipschitz and there exists a completely Lipschitz mapping $g_{l,q} : \Re^n \times C^0([-r,0];U) \to \Re^m$ and a constant $G \geq 0$ satisfying the following inequality for all $(x,u) \in \Re^n \times C^0([-r,0];U)$ :*

$$\left| g_{l,q}(x,u) \right| \leq G|x| + G\|u\|_r \tag{3.22}$$

*Finally, assume that for every $(x,u) \in \Re^n \times C^0([-r,+\infty);U)$ the solution $x(t)$ of (1.2) with initial condition $x(0) = x_0$ corresponding to $u \in C^0([-r,+\infty);U)$ satisfies $\frac{d}{dt}p_{l,q}(x(t),T_r(t)u) = g_{l,q}(x(t),T_r(t)u)$ for all $t \geq 0$ for which the solution exists.*

*If $l \geq 1$ is sufficiently large, then for every $\mu > 0$, there exists $\tilde{\sigma} \in KL$ such that for every $(x_0, w_0) \in \Re^n \times C^0([-r,0];\Re^m)$ the solution $(x(t), w(t))$ of (1.1) with*

$$u(t) = \Pr_U(w(t)) \tag{3.23}$$

$$\dot{w}(t) = g_{l,q}(x(t), T_r(t)u) - \mu(w(t) - p_{l,q}(x(t), T_r(t)u)) \tag{3.24}$$

*with initial condition $x(0) = x_0$, $T_r(0)w = w_0$ satisfies estimate (2.16), i.e., the dynamic feedback law (3.23), (3.24) achieves global stabilization of system (1.1). Moreover, if there exist constants $M, \omega > 0$ such that estimate (2.6') holds instead of (2.6) and $l \geq 1$ is sufficiently large then for every $\mu > 0$ there exist $\tilde{M}, \tilde{\omega} > 0$ such that for every $(x_0, w_0) \in \Re^n \times C^0([-r,0];\Re^m)$ the solution $(x(t), w(t))$ of (1.1), (3.23), (3.24) with initial condition $x(0) = x_0$, $T_r(0)w = w_0$ satisfies estimate (2.16) with $\tilde{\sigma}(s,t) := \tilde{M}\exp(-\tilde{\omega}t)s$, i.e., the dynamic feedback law (3.23), (3.24) achieves global exponential stabilization of system (1.1).*

**Remark 3.5:** It should be emphasized that the hypotheses of Corollary 3.4 usually hold if $k \in C^2(\Re^n; U)$ and $f \in C^2(\Re^n \times U; \Re^n)$ with $\left|\frac{\partial f}{\partial x}(x,u)\right| + \left|\frac{\partial f}{\partial u}(x,u)\right| \leq K$.

**Proof of Corollary 3.4:** By virtue of the assumptions and inequalities (3.20), (3.21) all hypotheses (A1-3) hold. Particularly, hypothesis (A2) holds with $a_1 = (1+\varepsilon)RK\frac{(LT)^{l+1}}{1-LT}$, $a_2 = (1+\varepsilon^{-1})RK\frac{(LT)^{l+1}}{1-LT}$ and $T = \frac{r}{q}$ for every $\varepsilon > 0$. It follows that (2.13) holds provided that

$$(\gamma + 1 + \gamma R)RK\frac{(LT)^{l+1}}{1-LT} < 1 \tag{3.25}$$

Since $LT < 1$, the above inequality is satisfied for sufficiently large $l \geq 1$. The conclusion is a consequence of Theorem 2.2. ◁

When the computation of the mapping $g_{l,q} : \Re^n \times C^0([-r,0]; U) \to \Re^m$ is difficult (due to high complexity of the formulae), one can use the following corollary (which is based on Theorem 2.3).

**Corollary 3.6:** *Consider systems (1.1) and (1.2) with $U = \Re^m$ and assume that hypothesis (A2) with a linear feedback $u = k'x$, where $k \in \Re^{n \times m}$. Further assume that inequalities (3.1a,b), hold. Define $\Phi_{l,q}(x,u) := P_{l,q}^{\delta_{-r}u} x$, for positive integers $l, q$ with $Lr < q$ and assume that $\Phi_{l,q} : \Re^n \times C^0([-r,0]; U) \to \Re^n$ is completely Lipschitz.*

*Let $\mu > 0$ be given. If $l \geq 1$ is sufficiently large, then there exists $\tilde{\sigma} \in KL$ such that for every $(x_0, u_0) \in \Re^n \times C^0([-r,0]; \Re^m)$ the solution $(x(t), w(t))$ of (1.1) with*

$$\dot{u}(t) = k'f\left(\Phi_{l,q}(x(t), T_r(t)u), u(t)\right) - \mu(u(t) - k'\Phi_{l,q}(x(t), T_r(t)u)) \tag{3.26}$$

*with initial condition $x(0) = x_0$, $T_r(0)u = u_0$ satisfies estimate (2.16), i.e., the dynamic feedback law (3.26) achieves global stabilization of system (1.1). Moreover, if there exist constants $M, \omega > 0$ such that estimate (2.6') holds instead of (2.6) and $l \geq 1$ is sufficiently large then there exist $\tilde{M}, \tilde{\omega} > 0$ such that for every $(x_0, u_0) \in \Re^n \times C^0([-r,0]; \Re^m)$ the solution $(x(t), w(t))$ of (1.1), (3.26) with initial condition $x(0) = x_0$, $T_r(0)u = u_0$ satisfies estimate (2.16) with $\tilde{\sigma}(s,t) := \tilde{M}\exp(-\tilde{\omega}t)s$, i.e., the dynamic feedback law (3.26) achieves global exponential stabilization of system (1.1).*

**Proof:** By virtue of the assumptions and inequality (3.16) all hypotheses (A1), (A2) and (A4) hold with $g(x,u) := k'f\left(\Phi_{l,q}(x,u), u(0)\right)$ for sufficiently large $l \geq 1$. Particularly, hypothesis (A4) holds with $a_1 = (1+\varepsilon)KR\max\{1, L\}\frac{(LT)^{l+1}}{1-LT}$, $a_2 = (1+\varepsilon^{-1})RK\max\{1, L\}\frac{(LT)^{l+1}}{1-LT}$ and $T = \frac{r}{q}$ for every $\varepsilon > 0$, where $R := |k|$.
It follows that (2.42) holds provided that

$$(\gamma + 1 + \gamma R)RK\max\{1, L\}\left(1 + \frac{1}{\mu}\right)\frac{(LT)^{l+1}}{1-LT} < 1 \tag{3.27}$$

Since $LT < 1$, the above inequality is satisfied for sufficiently large $l \geq 1$. The conclusion is a consequence of Theorem 2.3. ◁

## 4. Illustrating Example

The following example illustrates the use of Corollary 3.4 and Corollary 3.6 for the scalar system

$$\dot{x}(t) = f(x(t)) + u(t-r)$$
$$x(t) \in \Re, u \in \Re \quad (4.1)$$

where $f \in C^2(\Re;\Re)$ is a globally Lipschitz function with $|f'(x)| \leq 1$ for all $x \in \Re$. The feedback law

$$k(x) := -(1+\kappa)x \quad (4.2)$$

where $\kappa > 0$ is a constant, satisfies hypothesis (A2) with $\gamma = \frac{1+\varepsilon}{\kappa}$, $\sigma(s,t) := s(1+\varepsilon^{-1})\exp\left(-\frac{\kappa}{2}t\right)$ and $R = 1+\kappa$, where $\varepsilon > 0$ is arbitrary. This may be shown directly by using the time derivative of the function $V(x) = \frac{1}{2}x^2$ and the growth condition $|f(x)| \leq |x|$ for all $x \in \Re$. Hypothesis (A2) can be used for the small-gain analysis presented in [11]. More specifically, inequality (2.6') holds with $\gamma = \frac{1+\varepsilon}{\kappa\sqrt{1-\varepsilon}}$, $M := 1+\varepsilon^{-1}$ and $\omega := \frac{\kappa\varepsilon}{2}$, where $\varepsilon \in (0,1)$ is arbitrary.

Utilizing the small-gain arguments in [11], it may be shown that the closed-loop system (4.1) with $u(t) = k(x(t))$ will be globally asymptotically stable provided that

$$r < \frac{\kappa}{(1+\kappa)(2+\kappa)} \quad (4.3)$$

We next assume that $r < 1$. Notice that inequalities (3.1a,b) hold with $L = 1$. Corollary 3.4 guarantees that the family of approximate predictors $p_{l,1}: \Re^n \times C^0([-r,0];\Re) \to \Re$ (parameterized by the integer $l \geq 1$) will result to stabilizing feedback laws provided that (3.25) holds, i.e., provided that

$$2(1+\kappa)^2 r^{l+1} < \kappa(1-r) \quad (4.4)$$

Indeed, the proof of Proposition 3.2 shows that the constant $K(q)$ involved in (3.16) and (3.20) satisfies $K = 1$ for $q = 1$. Notice that the hypotheses of Corollary 3.4 are satisfied

- for $l = 1$ with

$$p_{1,1}(x,u) = -(1+\kappa)\left(x + rf(x) + \int_0^r u(\tau-r)d\tau\right) \text{ and } g_{1,1}(x,u) = -(1+\kappa)\left(f(x) + rf'(x)f(x) + rf'(x)u(-r) + u(0)\right)$$

- for $l = 2$ with

$$p_{2,1}(x,u) = -(1+\kappa)\left(x + \int_0^r f\left(x + \tau f(x) + \int_0^\tau u(s-r)ds\right)d\tau + \int_0^r u(\tau-r)d\tau\right) \text{ and}$$

$$g_{2,1}(x,u) = -(1+\kappa)\left(f(x) + u(0) + \int_0^r f'\left(x + \tau f(x) + \int_0^\tau u(s-r)ds\right)(f(x) + \tau f'(x)f(x) + \tau f'(x)u(-r) + u(\tau-r))d\tau\right)$$

If $\kappa = 3$, then system (4.1) will be stabilized by the "no prediction" feedback $u(t) = k(x(t))$ for $r < \frac{3}{20}$. On the other hand, inequality (4.4) shows that for every $\mu > 0$ the dynamic feedback (3.23), (3.24) with $l = q = 1$ will achieve global exponential stabilization for $r < \frac{\sqrt{393}-3}{64} \approx 0.2628$, i.e., higher values for the delay than $\frac{3}{20}$ are allowed. Moreover, inequality (4.4) shows that for every $\mu > 0$ the dynamic feedback (3.23), (3.24) with $q = 1$ will achieve global exponential stabilization for $r < r_{\max}(l)$, where $r_{\max}(l) \to 1$ as $l \to +\infty$. For example, for $l = 2$ global exponential stabilization is achieved for $r < 0.386$.

Other values for $q$ have to be used for the case $r \geq 1$. However, the formulae become very complicated for high values for the integers $l, q$. For $l = 1$, $q = 2$ we obtain the formulae:

$$p_{1,2}(x,u) = -(1+\kappa)\left(x + \frac{r}{2}f(x) + \int_0^r u(\tau - r)d\tau + \frac{r}{2}f\left(x + \frac{r}{2}f(x) + \int_0^{r/2} u(\tau - r)d\tau\right)\right)$$

$$g_{1,2}(x,u) = -(1+\kappa)\left(f(x) + \frac{r}{2}f'(x)f(x) + \frac{r}{2}f'(x)u(-r) + u(0)\right)$$

$$-(1+\kappa)\frac{r}{2}f'\left(x + \frac{r}{2}f(x) + \int_0^{r/2} u(\tau - r)d\tau\right)\left(f(x) + u\left(-\frac{r}{2}\right) + \frac{r}{2}f'(x)f(x) + \frac{r}{2}f'(x)u(-r)\right)$$

Exact computation of the constant $K(q)$ involved in (3.16) and (3.20) in conjunction with (3.25) shows that for every $\mu > 0$ the dynamic feedback (3.23), (3.24) with $q = 2$ will achieve global exponential stabilization provided that

$$2(1+\kappa)^2 r^{l+1}\left[1 + \exp\left(\frac{\sqrt{2}+1}{4}r\right) + \exp\left(\frac{r}{2}\right) + \frac{1}{2^l}\frac{r^{l+1}}{2-r}\right] < 2^l \kappa(2-r) \qquad (4.5)$$

For the case $\kappa = 3$, $l = 1$, the above inequality holds if $r < 0.3058$. Again, inequality (4.5) shows that for every $\mu > 0$ the dynamic feedback (3.23), (3.24) with $q = 2$ will achieve global exponential stabilization for $r < r_{\max}(l)$, where $r_{\max}(l) \to 2$ as $l \to +\infty$.

Corollary 3.6 can be used as well. For the case $l = 2$, $q = 2$ we have:

$$\Phi_{2,2}(x,u) = x_1 + \int_0^{r/2} f\left(x_1 + \tau f(x_1) + \int_{r/2}^{\tau + r/2} u(s-r)ds\right)d\tau + \int_{r/2}^r u(\tau - r)d\tau$$

$$x_1 = x + \int_0^{r/2} f\left(x + \tau f(x) + \int_0^\tau u(s-r)ds\right)d\tau + \int_0^{r/2} u(\tau - r)d\tau$$

By virtue of Corollary 3.6 and inequality (3.27) the dynamic feedback

$$\dot{u}(t) = -(1+\kappa)\left(f(\Phi_{2,2}(x(t), T_r(t)u)) + u(t)\right) - \mu(u(t) + (1+\kappa)\Phi_{2,2}(x(t), T_r(t)u))$$

will achieve exponential stabilization provided that $r < 2$ and

$$(1+\kappa)^2\left[1 + \exp\left(\frac{\sqrt{2}+1}{4}r\right) + \exp\left(\frac{r}{2}\right) + \frac{1}{4}\frac{r^3}{2-r}\right]r^3 < \frac{2\mu}{1+\mu}\kappa(2-r) \qquad (4.6)$$

For the case $\kappa = 3$, $\mu = 100$ the above inequality holds if $r < 0.5284$. Therefore, the use of simple predictor formulae, allowed a 252% increase in the value of the maximum allowable delay compared to the use of the "no prediction" feedback $u(t) = k(x(t))$.

The results can be applied in a similar way, to triangular systems of the form (1.3), where $f_i \in C^2(\Re^i; \Re)$ ($i = 1,...,n$) are globally Lipschitz functions with $|\nabla f_i(x)| \leq L$ for all $x \in \Re^i$ ($i = 1,...,n$). For such systems there exists a linear feedback for which hypothesis (A2) and inequality (2.6') hold (see [12]).

## 5. Concluding Remarks

In this work, sufficient conditions for global stabilization of nonlinear systems with delayed input by means of "approximate predictors" are presented. The approximate predictor is a notion introduced in the present work and roughly speaking is a mapping which approximates the exact values of the stabilizing input for the corresponding system with no delay. A systematic procedure for the construction of families of approximate predictors is provided for globally Lipschitz systems: the construction is based on successive approximations on appropriate time intervals. The resulting stabilizing feedback for the system with delayed input can be implemented by means of a dynamic distributed delay feedback law. An illustrating example showed the efficiency of the proposed control strategy for various predictor schemes.

Future research will address the important open problem of applying numerical methods for the construction of approximate predictors. Indeed, the recent work presented in [1] can be an alternative way for constructing approximate predictors for nonlinear systems which are not necessarily globally Lipschitz.

## **Appendix**

**Proof of Proposition 3.2:** For notational convenience we set $a := \dfrac{(LT)^{l+1}}{1-LT}$. Define the sequence:

$$g_i := \left| \underbrace{Q_{T,u_i}^k \ldots Q_{T,u_1}^k}_{i\ times} x - \phi(iT, x; u)) \right| \tag{A1}$$

By virtue of inequality (3.8) and definition (3.15) this sequence satisfies:

$$g_1 \le a(|x|+\|u\|) \text{ and } g_m = \left|P_{l,q}^u x - \phi(r,x;u)\right| \tag{A2}$$

where $\|u\| := \sup_{0\le\tau\le r} |u(\tau)|$. (A2) shows that inequality (3.16) holds with $K=1$ for the case $q=1$.

Next assume that $q \ge 2$. Inequality (3.8) implies the following recursive relation for $i = 1,\ldots,q-1$:

$$g_{i+1} := \left|\underbrace{Q^l_{T,u_{i+1}} Q^l_{T,u_i} \ldots Q^l_{T,u_1}}_{i+1 \text{ times}} x - \phi((i+1)T,x;u))\right| =$$

$$= \left|\underbrace{Q^l_{T,u_{i+1}} Q^l_{T,u_i} \ldots Q^l_{T,u_1}}_{i+1 \text{ times}} x - \phi(T,\phi(iT,x;u);u_{i+1}))\right| \le a\left(\left|\underbrace{Q^l_{T,u} \ldots Q^l_{T,u}}_{i \text{ times}} x\right| + \|u\|\right) + \exp(LT)\, g_i \tag{A3}$$

$$\le (a+\exp(LT))g_i + a(|\phi(iT,x;u)| + \|u\|)$$

Inequality (3.1b) guarantees that inequality (2.5) holds with $L$ replaced by $\dfrac{1+\sqrt{2}}{2}L$. Consequently, we obtain from (A3) for $i = 1,\ldots,q-1$:

$$g_{i+1} \le (a+\exp(LT))g_i + a(\exp(ipLT)+1)(|x|+\|u\|) \tag{A4}$$

where $p := \dfrac{1+\sqrt{2}}{2}$. Using (A4) and (A2) we can obtain the following estimate for $m \ge 2$:

$$\left|P_{l,q}^u x - \phi(r,x;u)\right| \le \left[(a+\exp(LT))^{q-1} + (\exp(pLr)+1)\frac{(a+\exp(LT))^{q-1}-1}{a+\exp(LT)-1}\right] a(|x|+\|u\|) \tag{A5}$$

Since $a = \dfrac{(LT)^{l+1}}{1-LT} \le \dfrac{(LT)^2}{1-LT}$, we obtain from (A5):

$$\left|P_{l,q}^u x - \phi(r,x;u)\right| \le \left[(b+\exp(LT))^{q-1} + (\exp(pLr)+1)\frac{(b+\exp(LT))^{q-1}-1}{b+\exp(LT)-1}\right] a(|x|+\|u\|) \tag{A6}$$

where $b := \dfrac{(LT)^2}{1-LT}$. It follows that inequality (3.16) holds with $K := (b+\exp(LT))^{q-1} + (\exp(pLr)+1)\dfrac{(b+\exp(LT))^{q-1}-1}{b+\exp(LT)-1}$ for the case $q \ge 2$. The proof is complete. ◁